\documentclass[12pt]{article}

\usepackage{amsmath, amsfonts, amssymb}
\usepackage{latexsym}
\usepackage{graphicx}
\usepackage{color}

\newtheorem{thm}{Theorem}[section]

\newcommand{\C}{\mathbb{C}}
\newcommand{\D}{\mathbb{D}}
\newcommand{\T}{\mathbb{T}}
\newcommand{\N}{\mathbb{N}}

\newcommand{\iy}{\infty}

\newcommand{\al}{\alpha}

\newcommand{\Ga}{\Gamma}
\newcommand{\de}{\delta}

\newcommand{\la}{\lambda}

\newcommand{\ph}{\varphi}

\newcommand{\si}{\sigma}
\newcommand{\tht}{\theta}

\newcommand{\n}{\|}

\newcommand{\ti}{\widetilde}
\newcommand{\ov}{\overline}

\renewcommand{\thefootnote}{\fnsymbol{footnote}}

\begin{document}

\begin{center}
{\Large \bf Borodin--Okounkov and Szeg\H{o} for Toeplitz \\[0.5ex] operators on model spaces}

\vspace{7mm}
{\Large Albrecht B\"ottcher}
\end{center}

\bigskip
\begin{quote}
\footnotesize{
We consider the determinants of compressions of Toeplitz operators to finite-dimensional model spaces
and establish analogues of the Borodin--Okounkov formula and the strong Szeg\H{o} limit theorem in this setting.}
\let\thefootnote\relax\footnote{\hspace*{-7.5mm} MSC 2010: 47B35, 30J10}
\let\thefootnote\relax\footnote{\hspace*{-7.5mm} Keywords: Toeplitz determinant, model space, Blaschke product, truncated Toeplitz operator}
\end{quote}

\section{Introduction and main results}

Although compressions of Toeplitz operators to model spaces have been studied for a long time, see, for example,
\cite{Nik}, \cite{Treil}, it was
Sarason's paper \cite{Sar} which initiated the recent increasing activity in research into
such operators\footnote[1]{These operators are
now called ``truncated Toeplitz operators'', although that name is already occupied
by the classical finite Toeplitz matrices. Moreover, I see a difference between truncation and compression.
However, since Donald Sarason is one of my mathematical
top heroes, I will not vote against that name. I will nevertheless not follow the custom
and will instead refer to these operators simply as Toeplitz operators on model spaces.}, see,
for instance, the survey \cite{GaRoss} and the ample list of references therein.
The number one theorem in classical Toeplitz matrices is Szeg\H{o}'s strong limit theorem,
and curiously, I have not seen the model space version of this theorem among the many results which
have so far been carried over from the classical setting to the model space level. In fact the
strong Szeg\H{o} limit theorem is a straightforward consequence of another great theorem,
namely, the Borodin--Okounkov formula. My favorite proof of the Borodin--Okounkov formula
is the one in \cite{Botok}, and the purpose of this note is to show that this proof works
equally well for Toeplitz operators on model spaces.

\medskip
Our context is the usual Hardy spaces of the unit disk $\D$ or, when interpreted as nontangential limits,
of the unit circle $\T$. We let $P$ stand for the orthogonal projection of $L^2$ onto $H^2$.
The Toeplitz operator $T(a)$ induced by a function $a \in L^\iy$ is the operator on $H^2$
which acts by the rule $T(a)f=P(af)$.
Let $u \in H^\iy$ be an inner function. The space $K_u:=H^2 \ominus uH^2$ is referred to as
the model space generated by $u$. We denote by $P_u$ and $Q_u=I-P_u$ the orthogonal
projections of $H^2$ onto $K_u$ and $uH^2$, respectively. It is well known
that $P_u=I-T(u)T(\ov{u})$, the bar denoting complex conjugation.
We are interested in the compression
of $T(a)$ to $K_u$, that is, in the operator $T_u(a)=P_u T(a)|K_u$.

\medskip
We will actually consider the matrix case. Thus, $a$ is supposed to be a matrix function in the
$\C^{m \times m}$-valued $L^\iy$, and $T(a)$ and $T_u(a)$ act on the $\C^m$-valued $H^2$ and $K_u$,
respectively.
The inner function $u$ remains scalar-valued.

\medskip
We make the following assumptions on $a$. It is required that $a$ is in the intersection of the Wiener algebra $W$
and the Krein algebra $K_{2,2}^{1/2,1/2}$, that is,
the Fourier coefficients $a_n$ satisfy
$\sum_{n=-\iy}^\iy \n a_n\n+\sum_{n=-\iy}^\iy n \n a_n\n^2 < \iy$,
where $\n \cdot\n$ is any matrix norm on $\C^{m \times m}$. We furthermore assume that $a$ has right and left
canonical Wiener--Hopf factorizations $a=w_-w_+=v_+v_-$ in $W \cap K_{2,2}^{1/2,1/2}$. This means that $w_+, v_+,
\ov{w_-}, \ov{v_-}$ and their inverses belong to $W \cap K_{2,2}^{1/2,1/2}\cap H^\iy$. In the scalar case ($m=1$),
the existence of such factorizations is guaranteed if $a$ has no zeros on $\T$ and
vanishing winding number about the origin. Our assumptions imply in particular that $T(a)$, $T(a^{-1})$, and $T(\ti{a})$
are invertible on $H^2$. Here and in what follows, $\ti{a}$ results from $a$ be reversal of the Fourier coefficients,
$\ti{a}(t):=a(1/t)$ for $t \in \T$.

\medskip
The Hankel operator $H(a)$ generated by $a \in L^\iy$ is defined on the space $H^2$ by $H(a)f=P(a\cdot (I-P)Jf)$,
where $J$ is the flip operator, $(Jf)(t)=(1/t)f(1/t)$ for $t \in \T$.
Put $b=v_-w_+^{-1}$ and $c=w_-^{-1}v_+$. Then $b$ and $c$ are in the Krein algebra and hence the Hankel
operators $H(b)$ and $H(\ti{c})$ are Hilbert--Schmidt operators. This implies that $H(b)H(\ti{c})$ is
in the trace class. As $T(b)=T(v_-)T(w_+^{-1})$ and $T(c)=T(w_-^{-1})T(v_+)$ are invertible, so also
is $I-H(b)H(\ti{c})=T(b)T(c)$.

\medskip
For $\al \in \D$, we define the inner functions $\mu_\al$ and $B_\al$ by
\[\mu_\al(z)=\frac{z-\al}{1-\ov{\al}z}, \quad B_\al(z)=\frac{-\ov{\al}}{|\al|}\frac{z-\al}{1-\ov{\al}z}
\quad (z \in \D),\]
with the convention to put $B_0(z)=z$. The space $K_u$ is known to be finite-dimensional if and only if
$u$ is a finite Blaschke product, that is, if and only if there are $\al_1, \ldots, \al_N$ in $\D$ such that
$u=B_{\al_1}\cdots B_{\al_N}$. We let $\si(u)$ denote the numbers $\al_1, \ldots, \al_N$, repeated
according to the number of times they appear in $u=B_{\al_1}\cdots B_{\al_N}$. Finally, as usual,
the geometric mean of a (matrix) function $\ph$ on $\T$ is defined by
\[G(\ph)=\exp(\log \det \ph)_0:=\exp\left(\frac{1}{2\pi}\int_0^{2\pi}\log\det\ph(e^{i\tht})d\tht\right).\]
Here is the model space version of the Borodin--Okounkov formula.

\begin{thm} \label{Theo 1.1}
If $u=B_{\al_1}\cdots B_{\al_N}$ is a finite Blaschke product, then
\begin{equation}
\det T_u(a)=\left(\prod_{\al \in \si(u)} G(a \circ \mu_{-\al})\right)\frac{\det (I-Q_u H(b)H(\ti{c})Q_u)}{\det (I-H(b)H(\ti{c}))}.
\label{1.1}
\end{equation}
\end{thm}

An alternative expression for the product of the numbers $G(a\circ\mu_{-\al})$ is
\begin{equation}
\prod_{\al \in \si(u)} G(a \circ \mu_{-\al})=\prod_{\al \in \si(u)} \det v_+(\al)\det v_-(1/\ov{\al}). \label{1.2}
\end{equation}
For $u(z)=z^N$, the products (\ref{1.2}) become $G(a)^N$ and (\ref{1.1}) turns into the classical
Borodin--Okounkov formula, which was originally established in \cite{Borodok}, reformulated, extended to the block case,  and
equipped with two new proofs in \cite{BasWid}, and with still another proof in \cite{Botok}.
For positive functions $a$, the formula was even already in \cite{GerCase}, which, however,
was not known to the authors of \cite{BasWid}, \cite{Borodok}, \cite{Botok} at the time they wrote
their papers. Taking into account that $Q_u=T(u)T(\ov{u})$ for an arbitrary inner function, it is easy to see that
\[\det(I-Q_u H(b)H(\ti{c})Q_u)=\det(I-H(\ov{u}b)H(\ti{c}\ti{u}))\]
for every inner function $u$.

\medskip
Now suppose $\{\al_j\}_{j=1}^\iy$ is a sequence of points in $\D$. Put
\[u_N(z)=\prod_{j=1}^N B_{\al_j}(z).\]
The following is a model space version of the strong Szeg\H{o} limit theorem.

\begin{thm} \label{Theo 1.2}
If $\sum_{j=1}^\iy (1-|\al_j|) = \iy$, then $u_N(z) \to 0$ for $z \in \D$, $Q_{u_N} \to 0$ strongly and
\begin{equation}
\lim_{N \to \iy} \det T_{u_N}(a)\prod_{\al \in \si(u_N)} G(a \circ \mu_{-\al})^{-1}= \frac{1}{\det (I-H(b)H(\ti{c}))}. \label{1.3}
\end{equation}
If  $\sum_{j=1}^\iy (1-|\al_j|) < \iy$, then $u_N(z)$ converges to the infinite Blaschke product
\[B(z)=\prod_{j=1}^\iy B_{\al_j}(z)\] for $z \in \D$, $Q_{u_N} \to Q_B$ strongly, and
\begin{equation}
\lim_{N \to \iy} \det T_{u_N}(a)\prod_{\al \in \si(u_N)} G(a \circ \mu_{-\al})^{-1}= \frac{\det(I-Q_B H(b)H(\ti{c})Q_B)}{\det (I-H(b)H(\ti{c}))}. \label{1.4}
\end{equation}
\end{thm}
Again, in the case where $u_N(z)=z^N$, this theorem implies that
\[\lim_{N \to \iy} T_{z^N}(a) G(a)^{-N} = \frac{1}{\det (I-H(b)H(\ti{c}))},\]
which is the classical Szeg\H{o}--Widom limit theorem, established by Szeg\H{o} \cite{Sz}
in the scalar case ($m=1$) and by Widom \cite{Wid} in the block case ($m \ge 1$).
Note that for $m=1$ we have
\[\frac{1}{\det (I-H(b)H(\ti{c}))}=\exp \sum_{k=1}^\iy k (\log a)_k (\log a)_{-k},\]
and that for $m \ge 1$ we may also write
\[\frac{1}{\det (I-H(b)H(\ti{c}))}=\det T(a) T(a^{-1}).\]
We refer to the books \cite{BoSi} and \cite{Simon} for more on this topic, including the history. Incidentally,
sequences of Toeplitz operators $T_{u_N}(a)$  with $u_{N+1}$ divisible by $u_N$ and with
$P_{u_N}$ converging strongly to $I$ appeared already in Treil's paper \cite{Treil} (and his results are also quoted
on p. 394 of \cite{BoSi}).

\section{Proofs}

We first prove Theorem \ref{Theo 1.1} and formula (\ref{1.2}).
Let $u$ be a finite Blaschke product. As shown in \cite{Botok} (or see \cite[p. 552]{BoSi} or \cite{BoWi}),
Jacobi's formula for the minors of the inverse matrix can be extended to identity minus trace class operators:
\[\det P_u (I-L)^{-1} P_u = \frac{\det (I-Q_u L Q_u)}{\det (I-L)}\]
whenever $L$ is of trace class and $I-L$ is invertible. This formula with $L=H(b)H(\ti{c})$ will give
Theorem~\ref{Theo 1.1} provided we can prove that
\begin{equation}
\det P_u(I-H(b)H(\ti{c}))^{-1}P_u=\det T_u(a) \prod_{\al \in \si(u)} G(a \circ \mu_{-\al})^{-1}. \label{2.1}
\end{equation}
It is readily seen that if $\ph \in H^\iy$, then
\begin{equation}
P_u T(\ph)=P_u T(\ph) P_u, \quad T(\ov{\ph}) P_u = P_u T(\ov{\ph}) P_u.\label{2.1a}
\end{equation}
Consequently,
\begin{eqnarray*}
& & P_u (I-H(b)H(\ti{c}))^{-1}P_u=P_u T(c)^{-1}T(b)^{-1}P_u\\
& & = P_u T(v_+^{-1})T(w_-)T(w_+)T(v_-^{-1})P_u
=T_u(v_+^{-1}) T_u(a) T_u(v_-^{-1}).
\end{eqnarray*}
Taking determinants, we see that the left-hand side of (\ref{2.1}) equals
\[\det T_u(a) /( \det T_u(v_+) \det T_u(v_-)).\]
We are so left with proving that
\begin{eqnarray}
& & \det T_u(v_+)=\prod_{\al \in \si(u)}\det v_+(\al), \quad \det T_u(v_-)=\prod_{\al \in \si(u)} \det v_-(1/\ov{\al}), \label{2.2}\\
& & \prod_{\al \in \si(u)}\det v_+(\al)\det v_-(1/\ov{\al})= \prod_{\al \in \si(u)} G(a \circ \mu_{-\al}). \label{2.3}
\end{eqnarray}
The determinant is the product of the eigenvalues. A complex number $\la$ is an eigenvalue of $T_u(v_+)$
if and only if $T_u(v_+)-\la I=T_u(v_+-\la I)$ is not invertible. We may think of $T_u(v_+-\la I)$ as an $m \times m$
block matrix whose blocks $T_u(v_+^{jk}-\la \de_{jk})$ are generated by scalar-valued functions. By virtue of~(\ref{2.1a}),
the blocks commute pairwise, and hence $T_u(v_+-\la I)$ is not invertible if and only if the block determinant
$\det T_u(v_+-\la I)$ is not invertible. Again by~(\ref{2.1a}), $\det T_u(v_+-\la I)=T_u(\det(v_+-\la I))$.
But the operator $T_u(\det(v_+-\la I))$ is known to be not invertible if and only if $\det(v_+(\al)-\la I)=0$ for
some $\al \in \si(u)$; see~\cite[p. 66]{Nik} or~\cite[Theorem 15(ii)]{GaRoss}. Equivalently, $T_u(\det(v_+-\la I))$ is
not invertible if and only if $\la$ is an eigenvalue of $v_+(\al)$ for some $\al \in \si(u)$. Thus, the set of the eigenvalues
of $T_u(v_+)$ is the union of the sets of the eigenvalues of $v_+(\al)$ for $\al \in \si(u)$, multiplicities taken into account.
This proves the first formula in~(\ref{2.2}). The second now follows from the equalities
\[\det T_u(v_-)=\ov{\det T_u(v_-^*)}=\prod_{\al \in \si(u)}\ov{\det v_-^*(\al)}=\prod_{\al \in \si(u)}\det v_-(1/\ov{\al}).\]
Finally, we have
\begin{eqnarray*}
& & \prod_{\al \in \si(u)}\det v_+(\al) \det v_-(1/\ov{\al})=
\prod_{\al \in \si(u)}\det(v_+\circ \mu_{-\al})(0) \det(v_- \circ \mu_{-\al})(\iy)\\
& & = \exp  \sum_{\al \in \si(u)}\Big(\log\det(v_+\circ \mu_{-\al})(0)+ \log\det(v_- \circ \mu_{-\al})(\iy)\Big)\\
& & = \exp  \sum_{\al \in \si(u)}\Big([\log\det(v_+\circ \mu_{-\al})]_0+ [\log\det(v_- \circ \mu_{-\al})]_0\Big)\\
& & = \exp  \sum_{\al \in \si(u)}[\log \det (a\circ \mu_{-\al})]_0=\prod_{\al \in \si(u)} G(a \circ \mu_{-\al}),
\end{eqnarray*}
which gives (\ref{2.3}) and completes the proof of Theorem \ref{Theo 1.1} and formula (\ref{1.2}).

\medskip
Once Theorem \ref{Theo 1.1} is available, Theorem \ref{Theo 1.2} is no surprise. Indeed, the assertions
concerning the limit of $u_N(z)$ are well known, and the theorem on the lower limits of model spaces
on page 35 of \cite{Nik} implies that $P_{u_N}$ converges strongly to $I$ if $u_N(z) \to 0$
and to $P_B$ if $u_N(z) \to B(z)$. Formulas~(\ref{1.3}) and~(\ref{1.4}) then result from Theorem~\ref{Theo 1.1}
and the continuity of the determinant on $I$ minus the trace ideal.

\section{Three Examples}

As already said, for $u(z)=z^N$ the term (\ref{1.2}) is simply $G(a)^N$. For general inner functions $u$,
it is less harmless. It suffices to illustrate things
in the simple case where $v_+(z)=1-vz$ with $|v| <1$. We put
\[G_u(v)=\prod_{\al \in \si(u)} v_+(\al)=\prod_{\al \in \si(u)} (1-v\al).\]

\medskip
{\bf Example 1.} Let $\al_j=1-1/j^2$ and $u_N(z)=\prod_{j=1}^N B_{\al_j}(z)$. Then
\begin{eqnarray*}
\log G_{u_N}(v) & = & \sum_{j=1}^N \log (1-v\al_j)=  \sum_{j=1}^N \log\left(1-v+\frac{v}{j^2}\right)\\
& = & N \log (1-v)+\sum_{j=1}^N \log \left(1+\frac{v}{1-v}\,\frac{1}{j^2}\right)\\
& = &  N \log (1-v)+\sum_{j=1}^\iy \log \left(1+\frac{v}{1-v}\,\frac{1}{j^2}\right)+O\left(\frac{1}{N}\right)
\end{eqnarray*}
and hence
\begin{eqnarray*}
G_{u_N}(v) & = & (1-v)^N \prod_{j=1}^\iy\left(1+\frac{v}{1-v}\,\frac{1}{j^2}\right)\left(1+O\left(\frac{1}{N}\right)\right)\\
& = & (1-v)^N \:\frac{\sinh\left(\pi \sqrt{\frac{v}{1-v}}\right)}{\pi \sqrt{\frac{v}{1-v}}}\left(1+O\left(\frac{1}{N}\right)\right).
\end{eqnarray*}

{\bf Example 2.} Now take $\al_j=1-1/j$ and  $u_N(z)=\prod_{j=1}^N B_{\al_j}(z)$. This time, with $q:=v/(1-v)$,
\begin{eqnarray*}
\log G_{u_N}(v) & = & N \log (1-v)+\sum_{j=1}^N \log\left(1+\frac{q}{j}\right)\\
& = & N\log(1-v)+\sum_{j=1}^N \left(\log\left(1+\frac{q}{j}\right)-\frac{q}{j}\right)+\sum_{j=1}^N\frac{q}{j},
\end{eqnarray*}
and this equals
\[N\log(1-v)+\sum_{j=1}^\iy \left(\log\left(1+\frac{q}{j}\right)-\frac{q}{j}\right)+O\left(\frac{1}{N}\right)+
q\left(\log N + C + O\left(\frac{1}{N}\right)\right),\]
where $C=0.5772\ldots$ is Euler's constant. It follows that
\[G_{u_N}(v)=(1-v)^N \,N^q\,e^{qC}\prod_{j=1}^\iy \left(1+\frac{q}{j}\right)e^{-q/j}\left(1+O\left(\frac{1}{N}\right)\right),\]
and taking into account that
\[\prod_{j=1}^\iy \left(1+\frac{q}{j}\right)e^{-q/j}=\frac{e^{-qC}}{\Ga(q+1)},\]
we arrive at the formula
\[G_{u_N}(v)=\frac{(1-v)^N\, N^{v/(1-v)}}{\Ga\left(\frac{1}{1-v}\right)}\left(1+O\left(\frac{1}{N}\right)\right).\]

{\bf Example 3.} The previous two examples raise the question whether the limits of $G_{u_{N+1}}(v)/G_{u_N}(v)$ and
$G_{u_N}(v)^{1/N}$ always exist. Surprisingly, the answer is NO. Since $G_{u_{N+1}}(v)/G_{u_N}(v)=1-v\al_{N+1}$, this
is clear for the quotient. To give a counterexample for the root, we construct a sequence $\{u_N\}$ with a subsequence
$\{u_{N_i}\}$ such that $G_{u_{N_i}}(v)^{1/N_i}$ alternately assumes two different values.
We take $u_N(z)=\prod_{j=1}^N B_{\al_j}(z)$
where $\al_j =r_j z_j$, $r_j \in (0,1)$, $z_j \in \T$, and
$\sum_{j=1}^\iy (1-r_j) < \iy$. Then
\begin{eqnarray*}
G_{u_N}(v) & = & \prod_{j=1}^N (1-vr_jz_j)=\prod_{j=1}^N (1-vz_j+vz_j(1-r_j))\\
& = & \prod_{j=1}^N (1-vz_j)\prod_{j=1}^N \left(1+\frac{vz_j}{1-vz_j}(1-r_j)\right)\\
& = &  \prod_{j=1}^N (1-vz_j)\prod_{j=1}^\iy \left(1+\frac{vz_j}{1-vz_j}(1-r_j)\right)\:(1+o(1)),
\end{eqnarray*}
and it is sufficient to choose $\{z_j\}_{j=1}^\iy$ so that the limit of $\prod_{j=1}^N(1-vz_j)^{1/N}$
does not exist. We successively take $z_j=-1$ or $z_j=1$ and denote by $f(N)$ the number of choices
of $z_j=1$ after $N$ steps. Here $f: \N \to \N$ may be any function such that
\begin{equation}
f(N-1) \le f(N) \le f(N-1)+1 \quad \mbox{for} \quad  N \ge 2. \label{3.1}
\end{equation}
Then
\[\prod_{j=1}^N(1-vz_j)^{1/N}=(1-v)^{f(N)/N}(1+v)^{(N-f(N))/N}=(1+v)\left(\frac{1-v}{1+v}\right)^{f(N)/N},\]
and we are left with searching a function satisfying (\ref{3.1}) such that $f(N)/N$ has no limit as $N \to \iy$.
Such functions obviously exist: start with $f(1)=1$, leave $f(N)$ constant until $f(N)/N=1/4$, then increase $f(N)$ successively by $1$
until $f(N)/N=1/2$, after that leave $f(N)$ again constant to reach $f(N)/N=1/4$, then increase
$f(N)$ anew by ones until $f(N)/N=1/2$, etc.
Here is this function explicitly. Every natural number $N \ge 3$ may uniquely be written as $N=2\cdot 3^k+\ell$ with $k \ge 0$
and $1 \le \ell \le 4\cdot 3^k$. We put
\[f(2\cdot 3^k+\ell)=\left\{\begin{array}{lll} 3^k & \mbox{for} & 1 \le \ell \le 2\cdot 3^k,\\
\ell-3^k & \mbox{for} & 2\cdot 3^k \le \ell \le 4\cdot 3^k, \end{array}\right.\]
and we also define $f(1)=f(2)=1$. Thus, our choice for $z_1$ is $1$, the following three choices are $z_2=z_3=z_4=-1$,
the following two are $z_5=z_6=1$, the following six $z_j$ are $-1$, the next six $z_j$ are $1$, and so on.
It can be verified straightforwardly that $f$ satisfies~(\ref{3.1}), and since
$f(N)/N=1/2$ for $N=2\cdot 3^k$ and $f(N)/N=1/4$ for $N=4\cdot 3^k$, the limit of $f(N)/N$ does
not exist.

\bigskip
Albrecht B\"ottcher

Fakult\"at f\"ur Mathematik

TU Chemnitz

09107 Chemnitz

Germany

\bigskip
{\tt aboettch@mathematik.tu-chemnitz.de}


\begin{thebibliography}{99}

\bibitem{BasWid}
E. Basor and H. Widom,
{\em On a Toeplitz determinant identity of Borodin and Okounkov}.
Integral Equations Operator Theory 37 (2000), 397–-401.

\bibitem{Borodok}
A. Borodin and A. Okounkov,
{\em A Fredholm determinant formula for Toeplitz determinants}.
Integral Equations Operator Theory 37 (2000), 386–-396.

\bibitem{Botok}
A. B\"ottcher,
{\em On the determinant formulas by Borodin, Okounkov, Baik, Deift, and Rains}.
Oper. Theory Adv. Appl. 135 (2002), 91--99.

\bibitem{BoSi}
A. B\"ottcher and B. Silbermann,
{\em Analysis of Toeplitz Operators}. Second edition, Springer-Verlag, Berlin, 2006.

\bibitem{BoWi}
A. B\"ottcher and H. Widom,
{\em Szeg\"o via Jacobi}.
Linear Algebra Appl. 419 (2006), 656–-667.

\bibitem{GaRoss}
S. R. Garcia and W. T. Ross,
{\em Recent progress in truncated Toeplitz operators}.
Preprint, arXiv: 1108.1858v4 [math.CV] 22 Jun 2012.

\bibitem{GerCase}
J. S. Geronimo and K. M. Case,
{\em Scattering theory and polynomials orthogonal on the unit circle}.
J. Math. Phys. 20 (1979), 299–-310.

\bibitem{Nik}
N. K. Nikolski,
{\em Treatise on the Shift Operator}.
Springer-Verlag, Berlin, 1986.

\bibitem{Sar}
D. Sarason,
{\em Algebraic properties of truncated Toeplitz operators}.
Oper. Matrices 1 (2007), 491–-526.

\bibitem{Simon}
B. Simon
{\em Orthogonal Polynomial on the Unit Circle. I}.
Amer. Math. Soc., Providence, RI, 2005.

\bibitem{Sz}
G. Szeg\H{o},
{\em On certain Hermitian forms associated with the Fourier series of a positive function}.
Festschrift Marcel Riesz, 228--238, Lund, 1952.

\bibitem{Treil}
S. R. Treil,
{\em Invertibility of a Toeplitz operator does not imply its invertibility by the projection method}.
Soviet Math. Dokl. 35 (1987), 103–-107.

\bibitem{Wid}
H. Widom,
{\em Asymptotic behavior of block Toeplitz matrices and determinants. II}.
Advances in Math. 21 (1976),  1–-29.

\end{thebibliography}
\end{document}